\documentclass[english,secthm,seceqn]{elsart1p}
\usepackage[T1]{fontenc}
\usepackage[latin9]{inputenc}
\usepackage{amsmath}
\usepackage{amssymb}

\makeatletter

\newcommand{\lyxdeleted}[3]{}

\usepackage[verbose]{cite}
\usepackage{babel}

\usepackage{mathptmx}

\makeatother

\usepackage{babel}

\begin{document}
\begin{frontmatter}

\title{Global dynamics of a cell mediated immunity in viral infection models
with distributed delays}

\author[Y]{Yukihiko Nakata}

\ead{nakata@bcamath.org}

\address[Y]{Basque Center for Applied Mathematics, Bizkaia Technology Park,
Building 500 E-48160 Derio, Spain }
\begin{abstract}
In this paper, we investigate global dynamics for a system of delay
differential equations which describes a virus-immune interaction
in \textit{vivo}. The model has two distributed time delays describing
time needed for infection of cell and virus replication. Our model
admits three possible equilibria, an uninfected equilibrium and infected
equilibrium with or without immune response depending on the basic
reproduction number for viral infection $R_{0}$ and for CTL response
$R_{1}$ such that $R_{1}<R_{0}$. It is shown that there always exists
one equilibrium which is globally asymptotically stable by employing
the method of Lyapunov functional. More specifically, the uninfected
equilibrium is globally asymptotically stable if $R_{0}\leq1$, an
infected equilibrium without immune response is globally asymptotically
stable if $R_{1}\leq1<R_{0}$ and an infected equilibrium with immune
response is globally asymptotically stable if $R_{1}>1$. The immune
activation has a positive role in the reduction of the infection cells
and the increasing of the uninfected cells if $R_{1}>1$.\end{abstract}
\end{frontmatter}

\section{Introduction}

The mathematical models, based on biological interactions, present
a framework which can be used to obtain new insights and to interpret
experimental data. Many authors have formulated mathematical models
which describe the dynamics of virus population in \textit{vivo} and
these provide advances in our understanding of HIV-1 (human immunodeficiency
virus 1) and other viruses, such as HBV (hepatitis B virus) and HCV
(hepatitis C virus) (see \cite{MR2255781,nowak1996population,MR1908737,herz1996viral,Huang2009,HZhu20092,Lv2009672,Mittler1998143,MR2340347,ChristinaBartholdy11152000,HZhu2009,MR2067116,MR2073965,MR2206233,MR2460257,Nelson2000201,Wang2007197,Wang2007532,Wodarz2002194,Gomez-Acevedo2010,Michael1,Perelson2002}
and the references therein). 

During viral infections, the host immune system reacts with antigen-specific
immune response. In particular, cytotoxic T lymphocytes (CTLs) play
a critical role in antiviral defense by attacking infected cells.
To investigate the relation between antiviral immune response and
virus load, Nowak and Bangham\cite{nowak1996population} developed
the following mathematical model. \begin{equation}
\begin{cases}
\frac{d}{dt}x(t)=s-dx(t)-kx(t)v(t),\\
\frac{d}{dt}y(t)=kx(t)v(t)-\delta y(t)-py(t)z(t),\\
\frac{d}{dt}v(t)=N\delta y(t)-\mu v(t),\\
\frac{d}{dt}z(t)=qy(t)z(t)-bz(t),\end{cases}\label{eq:standard}\end{equation}
where $x(t)$ denotes the concentration of uninfected target cells
at time $t$, $y(t)$ denotes the concentration of infected cells
that produce virus at time $t$, $v(t)$ denotes the concentration
of virus at time $t$ and $z(t)$ denotes the abundance of virus-specific
CTLs. Uninfected cells are produced at a constant rate $s$ and die
at rate $dx(t)$. Infected cells are produced from uninfected cells
and virus at rate $kx(t)v(t)$ and die at rate $\delta y(t)$. Free
virus is produced from uninfected cells at rate $N\delta y(t)$, where
$N$ denotes the total number of virus particles from one cell, and
die at rate $\mu v(t)$. The rate of CTL proliferation is given by
$qy(t)z(t)$ and decay at rate $bz(t)$ in the absence of stimulation
by the infected cells. Infected cells are killed by CTLs at rate $py(t)z(t)$.
All parameters are positive constant.

Korobeinikov \cite{MR2255781} studied global properties of a basic
viral infection model which ignores immunity ((\ref{eq:standard})
with $p=0$). By assuming that the incidence rate of infection is
given by a functional form, more general viral infection models are
proposed and investigated (see \cite{MR2340347,Huang2009}). Wodarz
et al. \cite{Wodarz2002194} considered a mathematical model for two
types immune responses. Murase et al. \cite{MR2206233} and Kajiwara
et al. \cite{MR2073965} studied stability of some mathematical models
for virus-immune interaction dynamics. Recently, Pr$\ddot{\textrm{u}}$ss
et al., \cite{MR2460257} showed that (\ref{eq:standard}) always
admits an equilibrium which is globally asymptotically stable by constructing
Lyapunov functions.

On the other hand, in modeling of many biological processes, time
delays are usually introduced for the purpose of accurate representations
of the phenomena. In virus dynamics, it has been assumed that new
virus particles are produced after the initial infection with a time
interval and this leads mathematical models by delay differential
equations. The estimated values of kinetic parameters are usually
changed by these delay differential equations (see \cite{herz1996viral,Mittler1998143,MR1908737,Nelson2000201}
and references therein). Mathematical analysis for these models is
necessary to obtain an integrated view for the virus dynamics in \textit{vivo}.
In particular, the global stability of a steady state for these models
will give us a detailed information and enhances our understanding
about the virus dynamics. 

In this paper, we introduce distributed (continuous) time delays to
(\ref{eq:standard}) and study its global dynamics. Let $h_{1}$ and
$h_{2}$ be positive constants and $f_{1}(\tau):\left[0,h_{1}\right]\to\mathbb{R}_{+}$
and $f_{2}(\tau):\left[0,h_{2}\right]\to\mathbb{R}_{+}$ be integrable
functions with $\int_{0}^{h_{1}}f_{1}(\tau)d\tau=\int_{0}^{h_{2}}f_{2}(\tau)d\tau=1$.
As in Mittler et al., \cite{Mittler1998143} and Nelson et al., \cite{MR1908737},
we assume that the infected cells $y(t)$ appear after the initial
infection with a time period $\tau$ and $\tau$ is distributed according
to $f_{1}(\tau)$ over the interval $\left[0,h_{1}\right]$, where
$h_{1}$ is the limit superior of the infection delay. In addition,
we assume that a time is needed for the virus production after a virions
enter a cell (see also \cite{Lv2009672,HZhu20092}). Thus, we also
assume the production delay $\tau$, which is distributed according
to $f_{2}(\tau)$ over the interval $\left[0,h_{2}\right]$, where
$h_{2}$ is the limit superior of this delay. Then, we obtain the
following viral infection model. \begin{equation}
\begin{cases}
\frac{d}{dt}x(t)=s-dx(t)-kx(t)v(t),\\
\frac{d}{dt}y(t)=k_{d}\int_{0}^{h_{1}}f_{1}(\tau)x(t-\tau)v(t-\tau)d\tau-\delta y(t)-py(t)z(t),\\
\frac{d}{dt}v(t)=N_{d}\delta\int_{0}^{h_{2}}f_{2}(\tau)y(t-\tau)d\tau-\mu v(t),\\
\frac{d}{dt}z(t)=qy(t)z(t)-bz(t).\end{cases}\label{eq:Mother-1}\end{equation}
The infection rate $k_{d}$ satisfies $k_{d}\leq k$ and the total
number of virus particles from one cell $N_{d}$ satisfies $N_{d}\leq N$,
if we incorporate the probability of surviving of the infected cells
and virus particles between the time for infection and for virus production,
respectively.

Stability analysis for (\ref{eq:Mother-1}) with discrete intracellular
delay was carried out by Li and Shu \cite{Michael1} and Zhu and Zou
\cite{HZhu2009}. Recently, based on Li and Shu \cite{Michael1},
Li and Shu \cite{li:2434} has investigated a viral infection model
with a general target cell dynamics, a nonlinear incidence rate and
distributed delay. Li and Shu \cite{Michael1,li:2434} showed that
their model always admits an equilibrium which is globally asymptotically
stable and it is necessary to have a logistic mitosis term in the
target cell dynamics for generating a periodic solution. On the other
hand, Zhu and Zou \cite{HZhu2009} established global stability of
an uninfected equilibrium and obtained sufficient conditions for local
asymptotic stability of two infected equilibria. However, since Li
and Shu \cite{Michael1,li:2434} did not consider the immune response
to the viral infection and Zhu and Zou \cite{HZhu2009} did not address
the global stability of the two infected equilibria for their model,
the global dynamics of (\ref{eq:Mother-1}) is still unclear and,
hence, our main aim is to establish the complete global dynamics.
We show that (\ref{eq:Mother-1}) has three possible equilibria, an
uninfected equilibrium and infected equilibrium with or without immune
response and always admits one equilibrium which is always globally
asymptotically stable. Moreover, it is shown that if the immune response
is activated, then the infected equilibrium with immune response is
globally stable. This implies that the immune response has a positive
role in the reduction of the infected cells.

The paper is organized as follows. In Section 2, we show the positivity
and ultimately boundedness of the solutions for (\ref{eq:Mother-1})
under suitable initial conditions. Then, we introduce two important
parameters, the basic reproduction number for viral infection $R_{0}$
and for CTL response $R_{1}$, defined by (\ref{eq:brn}) and (\ref{eq:brn2}),
respectively, and three possible equilibria for (\ref{eq:Mother-1}).
In Section 3, we establish global asymptotic stability of these equilibria
by constructing Lyapunov functional. It is shown that (\ref{eq:Mother-1})
always admits one equilibrium which is globally asymptotically stable
and, hence, we obtain the complete global dynamics of (\ref{eq:Mother-1}).
In Section 4, we study discrete delay model which was considered by
Zhu and Zou \cite{HZhu2009} and show that two infected equilibria
of their model is not only locally asymptotically stable but also
globally asymptotically stable. In Section 5, we offer a brief discussion.

\section{Preliminary results}

\subsection{Positivity and boundedness of the solutions}

To investigate the dynamics of (\ref{eq:Mother-1}), we set a suitable
phase space. Let $\overline{h}=\max\left\{ h_{1},h_{2}\right\} $.
We denote by $C=C([-\overline{h},0],\mathbb{R})$ the Banach space
of continuous functions mapping the interval $[-\overline{h},0]$
into $\mathbb{R}$ equipped with the sup-norm. The nonnegative cone
of $C$ is defined as $C_{+}=C([-\overline{h},0],\mathbb{R}_{+})$.
From the biological meanings, the initial conditions for (\ref{eq:Mother-1})
are\begin{equation}
x(\theta)=\varphi_{1}(\theta),y(\theta)=\varphi_{2}(\theta),v(\theta)=\varphi_{3}(\theta),z(0)=z_{0}\text{ for }\theta\in[-\overline{h},0],\label{eq:initial}\end{equation}
 where $\varphi_{i}\in C_{+},i=1,2,3$ and $z_{0}\geq0$.
\begin{lem}
\label{lem:positivity-1}Every solution of (\ref{eq:Mother-1}) with
(\ref{eq:initial}) is nonnegative for $t>0$. Every solution of (\ref{eq:Mother-1})
with (\ref{eq:initial}) is positive for $t>\overline{h}$ if $z_{0}>0$
and either

i) $\varphi_{2}(0)+\int_{0}^{h_{1}}f_{1}(\tau)\varphi_{1}(-\tau)\varphi_{3}(-\tau)d\tau>0$,
or

ii) $\varphi_{3}(0)+\int_{0}^{h_{2}}f_{2}(\tau)\varphi_{2}(-\tau)d\tau>0$.

Furthermore, every solution is bounded above by some positive constant
for sufficiently large $t$.\end{lem}
\begin{pf}
The solution $(x(t),y(t),v(t),z(t))$ of (\ref{eq:Mother-1}) with
(\ref{eq:initial}) exists and is unique on its maximal interval of
existence $(0,\sigma)$ for some $\sigma>0$. We see that $x(t)>0$
for all $t\in(0,\sigma)$. Indeed, this follows from that $\frac{d}{dt}x(t)=s>0$
for any $t\in(0,\sigma)$ when $x(t)=0$ from the first equation of
(\ref{eq:Mother-1}). It also holds that \[
z(t)=z_{0}\textrm{e}^{\int_{0}^{t}\left(qy(s)-b\right)ds}\geq0,\]
if $z_{0}\geq0$. In particular, $z(t)>0$ if $z_{0}>0$.

Let us show the nonnegativity of $y(t)$ and $v(t)$. Since we have
\begin{equation}
\begin{cases}
y(t) & =\left(\varphi_{2}(0)+k_{d}\int_{0}^{t}\int_{0}^{h_{1}}f_{1}(\tau)x(s-\tau)v(s-\tau)d\tau\textrm{e}^{\int_{0}^{s}\left(\delta+pz(u)\right)du}ds\right)\textrm{e}^{-\int_{0}^{t}\left(\delta+pz(s)\right)ds},\\
v(t) & =\left(\varphi_{3}(0)+N_{d}\delta\int_{0}^{t}\int_{0}^{h_{2}}f_{2}(\tau)y(s-\tau)d\tau\textrm{e}^{\mu s}ds\right)\textrm{e}^{-\mu t},\end{cases}\label{eq:int_pos}\end{equation}
from (\ref{eq:Mother-1}), $y(t)\geq0$ and $v(t)\geq0$ for $t>0$.
Now we show $y(t)>0$ and $v(t)>0$ for $t>\overline{h}$, if i) or
ii) holds.

First, we assume that i) holds. Suppose that there exists a $t_{1}$
such that $y(t_{1})=0$. Then, from (\ref{eq:int_pos}), \[
\varphi_{2}(0)+k_{d}\int_{0}^{t_{1}}\int_{0}^{h_{1}}f_{1}(\tau)x(s-\tau)v(s-\tau)d\tau\textrm{e}^{\int_{0}^{s}\left(\delta+pz(u)\right)du}ds=0,\]
follows. This leads a contradiction to i). Thus, we obtain \begin{equation}
y(t)>0\text{ for }t>0.\label{eq:y1_positive}\end{equation}
Next, we suppose that there exists a $t_{2}>h_{2}$ such that $v(t_{2})=0$.
Then, from (\ref{eq:int_pos}), \[
\varphi_{3}(0)+N_{d}\delta\int_{0}^{t_{2}}\int_{0}^{h_{2}}f_{2}(\tau)y(s-\tau)d\tau\textrm{e}^{\mu s}ds=0,\]
follows. On the other hand, we have \[
\int_{0}^{t}\int_{0}^{h_{2}}f_{2}(\tau)y(s-\tau)d\tau ds>0,\text{ for }t>h_{2},\]
by (\ref{eq:y1_positive}). This gives a contradiction. Thus, $v(t)>0$
for $t>h_{2}$. Similarly, we see $v(t)>0$ for $t>0$ and $y(t)>0$
for $t>h_{1}$ if ii) holds.

Now we show the boundedness of each solution. Let\[
G(t)=\frac{k_{d}}{k}\int_{0}^{h_{1}}f_{1}(\tau)x(t-\tau)d\tau+y(t)+\frac{p}{q}z(t),\]
then we see \begin{align*}
\frac{d}{dt}G(t) & =\left(\frac{k_{d}}{k}\int_{0}^{h_{1}}f_{1}(\tau)\left(s-dx(t-\tau)-kx(t-\tau)v(t-\tau)\right)d\tau\right)\\
 & +\left(k_{d}\int_{0}^{h_{1}}f_{1}(\tau)x(t-\tau)v(t-\tau)d\tau-\delta y(t)-py(t)z(t)\right)+\left(py(t)z(t)-\frac{p}{q}bz(t)\right)\\
 & =\frac{sk_{d}}{k}-\frac{dk_{d}}{k}\int_{0}^{h_{1}}f_{1}(\tau)x(t-\tau)d\tau-\delta y(t)-\frac{p}{q}bz(t).\end{align*}
Therefore, it follows that\begin{align*}
\frac{d}{dt}G(t) & \leq\frac{sk_{d}}{k}-\min\left\{ d,\delta,b\right\} G(t),\end{align*}
which implies that $x(t),y(t)$ and $z(t)$ are uniformly bounded
on $(0,\sigma)$. Then, $v(t)$ is also uniformly bounded on $(0,\sigma)$.
Finally, it follows that $(x(t),y(t),v(t),z(t))$ exists and is unique
and positive for any $t>\overline{h}$.\qed\end{pf}
\begin{rem}
$y(t)$ and $v(t)$ of (\ref{eq:Mother-1}) with (\ref{eq:initial})
are identically zero for $t>0$, if \[
\varphi_{2}(0)=\varphi_{3}(0)=\int_{0}^{h_{1}}f_{1}(\tau)\varphi_{1}(-\tau)\varphi_{3}(-\tau)d\tau=\int_{0}^{h_{2}}f_{2}(\tau)\varphi_{2}(-\tau)d\tau=0.\]

\end{rem}

\subsection{Possible equilibria}

In this subsection, we show that (\ref{eq:Mother-1}) has three possible
equilibria. Existence of these equilibria is determined by a combination
of two threshold parameters\begin{equation}
R_{0}=\frac{s}{d\frac{\mu}{k_{d}N_{d}}},\label{eq:brn}\end{equation}
and \begin{equation}
R_{1}=\frac{s}{d\frac{\mu}{k_{d}N_{d}}+\frac{k}{k_{d}}\delta\frac{b}{q}}.\label{eq:brn2}\end{equation}
$R_{0}$ and $R_{1}$ are called the basic reproduction number for
viral infection and for CTL response, respectively (see Gomez-Acevedo
et al. \cite{Gomez-Acevedo2010}). In particular, $R_{0}$ denotes
the average number of secondary virus produced from a single virus
for (\ref{eq:Mother-1}).
\begin{thm}
\label{thm:DFE-1}For (\ref{eq:Mother-1}), there exist an uninfected
equilibrium \begin{equation}
E_{0}=(x_{0},0,0,0),x_{0}=\frac{s}{d},\label{eq:DFE}\end{equation}
an infected equilibrium without immune response \begin{equation}
E_{1}=(x_{1}^{*},y_{1}^{*},v_{1}^{*},0)=\left(\frac{\mu}{k_{d}N_{d}},\frac{k_{d}}{k\delta}\left(s-d\frac{\mu}{k_{d}N_{d}}\right),\frac{k_{d}N_{d}}{k\mu}\left(s-d\frac{\mu}{k_{d}N_{d}}\right),0\right),\label{eq:EE1}\end{equation}
if $R_{0}>1$, and an infected equilibrium with immune response \begin{equation}
E_{2}=(x_{2}^{*},y_{2}^{*},v_{2}^{*},z_{2}^{*})=\left(\frac{s}{d+k\frac{N_{d}\delta}{\mu}\frac{b}{q}},\frac{b}{q},\frac{N_{d}\delta}{\mu}\frac{b}{q},\frac{\delta}{p}\left(\frac{k_{d}N_{d}x_{2}^{*}}{\mu}-1\right)\right),\label{eq:EE2}\end{equation}
if $R_{1}>1$.\end{thm}
\begin{pf}
First of all, we see that (\ref{eq:Mother-1}) always has the uninfected
equilibrium $E_{0}$. To find other equilibria, we consider the following
equations\begin{equation}
\begin{cases}
0=s-dx^{*}-kx^{*}v^{*},\\
0=k_{d}x^{*}v^{*}-\delta y^{*}-py^{*}z^{*},\\
0=N_{d}\delta y^{*}-\mu v^{*},\\
0=qy^{*}z^{*}-bz^{*}.\end{cases}\label{eq:steady}\end{equation}
Assume that there exists an equilibrium $E_{1}=(x_{1}^{*},y_{1}^{*},v_{1}^{*},0)$
with $x_{1}^{*}>0,y_{1}^{*}>0,v_{1}^{*}>0$. From the third equation
of (\ref{eq:steady}), we see \begin{equation}
y_{1}^{*}=\frac{\mu}{N_{d}\delta}v_{1}^{*}.\label{eq:yv}\end{equation}
Substituting (\ref{eq:yv}) into the second equation of (\ref{eq:steady})
gives\[
0=k_{d}x_{1}^{*}v_{1}^{*}-\frac{\mu}{N_{d}}v_{1}^{*}.\]
Hence,\[
x_{1}^{*}=\frac{\mu}{k_{d}N_{d}}.\]
Then from the first equation of (\ref{eq:steady}), it follows\[
v_{1}^{*}=\frac{s-dx_{1}^{*}}{kx_{1}^{*}}=\frac{k_{d}N_{d}}{k\mu}\left(s-d\frac{\mu}{k_{d}N_{d}}\right).\]
$v_{1}^{*}$ and $y_{1}^{*}$ is positive, if $R_{0}>1$. Consequently,
there exists the infected equilibrium $E_{1}$ if $R_{0}>1$. 

Next, we assume that there exists an equilibrium $E_{2}=(x_{2}^{*},y_{2}^{*},v_{2}^{*},z_{2}^{*})$
with $x_{2}^{*}>0,y_{2}^{*}>0,v_{2}^{*}>0,z_{2}^{*}>0$. We have\begin{equation}
y_{2}^{*}=\frac{b}{q},v_{2}^{*}=\frac{N_{d}\delta}{\mu}y_{2}^{*},\label{eq:yv2}\end{equation}
from the forth and third equations of (\ref{eq:steady}), respectively.
Then, we have \begin{equation}
x_{2}^{*}=\frac{s}{d+kv_{2}^{*}}=\frac{s}{d+k\frac{N_{d}\delta}{\mu}\frac{b}{q}},\label{eq:x2-1}\end{equation}
and \[
z_{2}^{*}=\frac{k_{d}x_{2}^{*}v_{2}^{*}-\delta y_{2}^{*}}{py_{2}^{*}},\]
from the first and second equations of (\ref{eq:steady}), respectively.
By (\ref{eq:yv2}) and (\ref{eq:x2-1}), we see \[
z_{2}^{*}=\frac{k_{d}x_{2}^{*}}{p}\frac{N_{d}\delta}{\mu}-\frac{\delta}{p}=\frac{\delta}{p}\left(\frac{k_{d}N_{d}x_{2}^{*}}{\mu}-1\right)=\frac{\delta}{p}\left(R_{1}-1\right).\]
Thus, $z_{2}^{*}$ is positive if $R_{1}>1$ and, hence, there exists
the infected equilibrium with immune response $E_{2}$. Consequently,
the proof is complete. \qed\end{pf}
\begin{rem}
\label{rem:equ_cond}For $R_{1}>1$, there exist three equilibria,
$E_{0}$, $E_{1}$ and $E_{2}$. Moreover, we have $x_{2}^{*}>x_{1}^{*}$
and \textup{$y_{1}^{*}>y_{2}^{*}$,} since\[
x_{2}^{*}-x_{1}^{*}=\frac{\mu}{k_{d}N_{d}}\left(R_{1}-1\right)>0,\]
and\[
y_{1}^{*}-y_{2}^{*}=\frac{k_{d}}{k\delta}\left(s-d\frac{\mu}{k_{d}N_{d}}-\frac{k\delta}{k_{d}}\frac{b}{q}\right)=\frac{k_{d}}{k\delta}\left(d\frac{\mu}{k_{d}N_{d}}+\frac{k\delta}{k_{d}}\frac{b}{q}\right)\left(R_{1}-1\right)>0,\]
follows. Therefore, for the equilibrium condition, the immune activation
has a positive role in the increasing of the uninfected cells and
the reduction of the infected cells.
\end{rem}
Thus, there exist three possible equilibria depending on the values
of $R_{0}$ and $R_{1}$ defined by (\ref{eq:brn}) and (\ref{eq:brn2}),
respectively. We see that $R_{0}>R_{1}$ always holds.

\section{Global asymptotic stability of three equilibria}

In this section, we study the global dynamics of (\ref{eq:Mother-1})
by employing the method of Lyapunov functional. Lyapunov functionals,
we construct here, are inspired by McClusky \cite{McCluskey201055}
for SIR epidemic models with distributed delay. From the following
result, we see that (\ref{eq:Mother-1}) always admits one equilibrium
which is globally asymptotically stable and hence, the global dynamics
of (\ref{eq:Mother-1}) is fully determined by $R_{0}$ and $R_{1}$.
\begin{thm}
\label{thm:Main}i) If $R_{0}\leq1$, then the uninfected equilibrium
$E_{0}$ for (\ref{eq:Mother-1}) is globally asymptotically stable.

ii) Assume that either i) or ii) in Lemma \ref{lem:positivity-1}
holds. If $R_{1}\leq1<R_{0}$, then the infected equilibrium without
immune response $E_{1}$ for (\ref{eq:Mother-1}) is globally asymptotically
stable.

iii) Assume that $z_{0}>0$ and either i) or ii) in Lemma \ref{lem:positivity-1}
holds. If $R_{1}>1$, then the infected equilibrium with immune response
$E_{2}$ for (\ref{eq:Mother-1}) is globally asymptotically stable.
\end{thm}
Before giving the proof of Theorem \ref{thm:Main}, we introduce some
notations. In the Lyapunov functionals, the following function is
useful. \[
g(x)=x-1-\ln x,\text{ for }x\in(0,+\infty).\]
 $g(x)$ has the global minimum at $x=1$ and $g(1)=0$.

For simplicity, we will use the following notation in the proof \begin{align*}
\tilde{x}_{t} & =\frac{x(t)}{x_{1}^{*}},\tilde{x}_{t,\tau}=\frac{x(t-\tau)}{x_{1}^{*}},\tilde{y}_{t}=\frac{y(t)}{y_{1}^{*}},\tilde{y}_{t,\tau}=\frac{y(t-\tau)}{y_{1}^{*}},\tilde{v}_{t}=\frac{v(t)}{v_{1}^{*}},\tilde{v}_{t,\tau}=\frac{v(t-\tau)}{v_{1}^{*}},\tilde{z}_{t}=\frac{z(t)}{z_{1}^{*}},\\
\overline{x}_{t} & =\frac{x(t)}{x_{2}^{*}},\overline{x}_{t,\tau}=\frac{x(t-\tau)}{x_{2}^{*}},\overline{y}_{t}=\frac{y(t)}{y_{2}^{*}},\overline{y}_{t,\tau}=\frac{y(t-\tau)}{y_{2}^{*}},\overline{v}_{t}=\frac{v(t)}{v_{2}^{*}},\overline{v}_{t,\tau}=\frac{v(t-\tau)}{v_{2}^{*}},\overline{z}_{t}=\frac{z(t)}{z_{2}^{*}},\end{align*}
for $\tau\in[0,\overline{h}]$.
\begin{pf}
i) We construct the following Lyapunov functional\begin{equation}
U_{0}(t)=\frac{k_{d}}{k}x_{0}g\left(\frac{x(t)}{x_{0}}\right)+y(t)+\frac{1}{N_{d}}v(t)+\frac{p}{q}z(t)+\overline{U}_{0}(t),\label{eq:Lyapunov0}\end{equation}
where\begin{align*}
\overline{U}_{0}(t) & =k_{d}\int_{0}^{h_{1}}f_{1}(\tau)\int_{t-\tau}^{t}x(s)v(s)dsd\tau+\delta\int_{0}^{h_{2}}f_{2}(\tau)\int_{t-\tau}^{t}y(s)dsd\tau.\end{align*}
We calculate the time derivative of $U_{0}(t)$ along the solutions
of (\ref{eq:Mother-1}). We see \begin{align}
\frac{d}{dt}\left[x_{0}g\left(\frac{x(t)}{x_{0}}\right)\right] & =x_{0}\left(\frac{1}{x_{0}}-\frac{1}{x(t)}\right)\left(s-dx(t)-kx(t)v(t)\right)\nonumber \\
 & =\left(1-\frac{x_{0}}{x(t)}\right)\left(dx_{0}-dx(t)-kx(t)v(t)\right)\nonumber \\
 & =\left(1-\frac{x_{0}}{x(t)}\right)dx(t)\left(\frac{x_{0}}{x(t)}-1\right)-\left(1-\frac{x_{0}}{x(t)}\right)kx(t)v(t)\nonumber \\
 & =-dx(t)\left(1-\frac{x_{0}}{x(t)}\right)^{2}-kx(t)v(t)+kx_{0}v(t).\label{eq:x1}\end{align}
Next, we obtain\begin{align}
 & \frac{d}{dt}\left(y(t)+\frac{1}{N_{d}}v(t)+\frac{p}{q}z(t)\right)\nonumber \\
= & k_{d}\int_{0}^{h_{1}}f_{1}(\tau)x(t-\tau)v(t-\tau)d\tau-\delta y(t)-py(t)z(t)\nonumber \\
 & +\frac{1}{N_{d}}\left(N_{d}\delta\int_{0}^{h_{2}}f_{2}(\tau)y(t-\tau)d\tau-\mu v(t)\right)+\frac{p}{q}\left(qy(t)z(t)-bz(t)\right)\nonumber \\
= & k_{d}\int_{0}^{h_{1}}f_{1}(\tau)x(t-\tau)v(t-\tau)d\tau-\delta y(t)+\delta\int_{0}^{h_{2}}f_{2}(\tau)y(t-\tau)d\tau-\frac{\mu}{N_{d}}v(t)-\frac{p}{q}bz(t).\label{eq:y1}\end{align}
Finally, we obtain \begin{align}
\frac{d}{dt}\overline{U}_{0}(t) & =k_{d}\int_{0}^{h_{1}}f_{1}(\tau)\left(x(t)v(t)-x(t-\tau)v(t-\tau)\right)d\tau+\delta\int_{0}^{h_{2}}f_{2}(\tau)\left(y(t)-y(t-\tau)\right)dsd\tau\nonumber \\
 & =k_{d}\left(x(t)v(t)-\int_{0}^{h_{1}}f_{1}(\tau)x(t-\tau)v(t-\tau)d\tau\right)+\delta\left(y(t)-\int_{0}^{h_{2}}f_{2}(\tau)y(t-\tau)d\tau\right).\label{eq:ca1}\end{align}

Consequently, by adding (\ref{eq:x1}), (\ref{eq:y1}) and (\ref{eq:ca1}),
we obtain \begin{align*}
\frac{d}{dt}U_{0}(t) & =-\frac{k_{d}d}{k}x(t)\left(1-\frac{x_{0}}{x(t)}\right)^{2}+k_{d}x_{0}v(t)-\frac{\mu}{N_{d}}v(t)-\frac{p}{q}bz(t)\\
 & =-\frac{k_{d}d}{k}x(t)\left(1-\frac{x_{0}}{x(t)}\right)^{2}+\left(k_{d}x_{0}-\frac{\mu}{N_{d}}\right)v(t)-\frac{p}{q}bz(t)\\
 & =-\frac{k_{d}d}{k}x(t)\left(1-\frac{x_{0}}{x(t)}\right)^{2}+\frac{\mu}{N_{d}}\left(R_{0}-1\right)v(t)-\frac{p}{q}bz(t)\leq0,\text{ for }R_{0}\leq1.\end{align*}
Hence, every solution of (\ref{eq:Mother-1}) tends to $M_{0}$, where
$M_{0}$ is the largest invariant subset in $\left\{ \frac{dU_{0}(t)}{dt}=0\right\} $
with respect to (\ref{eq:Mother-1}). We show that $M_{0}$ consists
of only the equilibrium $E_{0}$. Let $\left(x(t),y(t),v(t),z(t)\right)$
be the solution with initial function in $M_{0}$. Then, from the
invariance of $M_{0}$, $x(t)=x_{0}$ and $z(t)=0$ for any $t$.
Now we have $\frac{d}{dt}x(t)=0$ and hence, it follows $v(t)=0$
for any $t$, from the first equation of (\ref{eq:Mother-1}). Then,
from the second equation of (\ref{eq:Mother-1}), we obtain $\lim_{t\to+\infty}y(t)=0$.
Therefore, the uninfected equilibrium $E_{0}$ is globally attractive.
Since we have $\frac{dU_{0}(t)}{dt}\leq0$ for $R_{0}\leq1$ and $U_{0}(t)\geq U_{0}(t)-\overline{U}_{0}(t)$,
the uninfected equilibrium $E_{0}$ is stable by Hale and Lunel \cite[Section 5, Corollary 3.1]{MR1243878}.
Hence, the uninfected equilibrium $E_{0}$ is globally asymptotically
stable for $R_{0}\leq1$.

ii) We construct the following Lyapunov functional \begin{equation}
U_{1}(t)=\frac{1}{kv_{1}^{*}}g\left(\frac{x(t)}{x_{1}^{*}}\right)+\frac{y_{1}^{*}}{k_{d}x_{1}^{*}v_{1}^{*}}g\left(\frac{y(t)}{y_{1}^{*}}\right)+\frac{v_{1}^{*}}{N_{d}\delta y_{1}^{*}}g\left(\frac{v(t)}{v_{1}^{*}}\right)+\frac{p}{k_{d}x_{1}^{*}v_{1}^{*}q}z(t)+\overline{U}_{1}(t),\label{Lyapunov}\end{equation}
 where\[
\overline{U}_{1}(t)=\int_{0}^{h_{1}}f_{1}(\tau)\int_{t-\tau}^{t}g\left(\frac{x(s)v(s)}{x_{1}^{*}v_{1}^{*}}\right)dsd\tau+\int_{0}^{h_{2}}f_{2}(\tau)\int_{t-\tau}^{t}g\left(\frac{y_{1}(s)}{y_{1}^{*}}\right)dsd\tau.\]
We calculate the time derivative of $U_{1}(t)$ along the positive
solutions of (\ref{eq:Mother-1}) and show that $\frac{dU_{1}(t)}{dt}\leq0$.
First, we have \begin{align*}
\frac{d}{dt}\left[g\left(\frac{x(t)}{x_{1}^{*}}\right)\right] & =\frac{1}{x_{1}^{*}}\left(1-\frac{x_{1}^{*}}{x(t)}\right)\left(s-dx(t)-kx(t)v(t)\right).\end{align*}
Since $s=dx_{1}^{*}+kx_{1}^{*}v_{1}^{*}$ holds, it follows\begin{align}
\frac{d}{dt}\left[g\left(\frac{x(t)}{x_{1}^{*}}\right)\right] & =\frac{1}{x_{1}^{*}}\left(1-\frac{x_{1}^{*}}{x(t)}\right)\left(dx_{1}^{*}+kx_{1}^{*}v_{1}^{*}-dx(t)-kx(t)v(t)\right)\nonumber \\
 & =\frac{1}{x_{1}^{*}}\left(1-\frac{x_{1}^{*}}{x(t)}\right)\left(dx_{1}^{*}-dx(t)+kx_{1}^{*}v_{1}^{*}-kx(t)v(t)\right)\nonumber \\
 & =-\frac{dx(t)}{x_{1}^{*}}\left(1-\frac{x_{1}^{*}}{x(t)}\right)^{2}+kv_{1}^{*}\left(1-\frac{x_{1}^{*}}{x(t)}\right)\left(1-\frac{x(t)v(t)}{x_{1}^{*}v_{1}^{*}}\right)\nonumber \\
 & =-\frac{dx(t)}{x_{1}^{*}}\left(1-\frac{x_{1}^{*}}{x(t)}\right)^{2}+kv_{1}^{*}\left(1-\frac{1}{\tilde{x}_{t}}\right)\left(1-\tilde{x}_{t}\tilde{v}_{t}\right)\nonumber \\
 & =-\frac{dx(t)}{x_{1}^{*}}\left(1-\frac{x_{1}^{*}}{x(t)}\right)^{2}+kv_{1}^{*}\left(1-\tilde{x}_{t}\tilde{v}_{t}-\frac{1}{\tilde{x}_{t}}+\tilde{v}_{t}\right).\label{eq:x2}\end{align}
Secondly, we compute \begin{align*}
\frac{d}{dt}\left[g\left(\frac{y(t)}{y_{1}^{*}}\right)\right] & =\frac{1}{y_{1}^{*}}\left(1-\frac{y_{1}^{*}}{y(t)}\right)\left(k_{d}\int_{0}^{h_{1}}f_{1}(\tau)x(t-\tau)v(t-\tau)d\tau-\delta y(t)-py(t)z(t)\right)\\
 & =\frac{1}{y_{1}^{*}}\left(1-\frac{y_{1}^{*}}{y(t)}\right)\left(\int_{0}^{h_{1}}f_{1}(\tau)\left(k_{d}x(t-\tau)v(t-\tau)-\delta y(t)\right)d\tau-py(t)z(t)\right).\end{align*}
Since we have $\delta=\frac{k_{d}x_{1}^{*}v_{1}^{*}}{y_{1}^{*}}$,
it follows \begin{align}
 & \frac{d}{dt}\left[g\left(\frac{y(t)}{y_{1}^{*}}\right)\right]\nonumber \\
 & =\frac{1}{y_{1}^{*}}\left(1-\frac{y_{1}^{*}}{y(t)}\right)\int_{0}^{h_{1}}f_{1}(\tau)\left(k_{d}x(t-\tau)v(t-\tau)-k_{d}x_{1}^{*}v_{1}^{*}\frac{y(t)}{y_{1}^{*}}\right)d\tau-\frac{1}{y_{1}^{*}}\left(1-\frac{y_{1}^{*}}{y(t)}\right)py(t)z(t)\nonumber \\
 & =\frac{k_{d}x_{1}^{*}v_{1}^{*}}{y_{1}^{*}}\left(1-\frac{y_{1}^{*}}{y(t)}\right)\int_{0}^{h_{1}}f_{1}(\tau)\left(\frac{x(t-\tau)v(t-\tau)}{x_{1}^{*}v_{1}^{*}}-\frac{y(t)}{y_{1}^{*}}\right)d\tau-\frac{1}{y_{1}^{*}}\left(1-\frac{y_{1}^{*}}{y(t)}\right)py(t)z(t)\nonumber \\
 & =\frac{k_{d}x_{1}^{*}v_{1}^{*}}{y_{1}^{*}}\left(1-\frac{1}{\tilde{y}_{t}}\right)\int_{0}^{h_{1}}f_{1}(\tau)\left(\tilde{x}_{t,\tau}\tilde{v}_{t,\tau}-\tilde{y}_{t}\right)d\tau-\frac{1}{y_{1}^{*}}\left(1-\frac{y_{1}^{*}}{y(t)}\right)py(t)z(t)\nonumber \\
 & =\frac{k_{d}x_{1}^{*}v_{1}^{*}}{y_{1}^{*}}\int_{0}^{h_{1}}f_{1}(\tau)\left(\tilde{x}_{t,\tau}\tilde{v}_{t,\tau}-\frac{\tilde{x}_{t,\tau}\tilde{v}_{t,\tau}}{\tilde{y}_{t}}-\tilde{y}_{t}+1\right)d\tau-\frac{1}{y_{1}^{*}}\left(py(t)z(t)-py_{1}^{*}z(t)\right).\label{eq:y2}\end{align}
Let us calculate the following \begin{eqnarray*}
\frac{d}{dt}\left[g\left(\frac{v(t)}{v_{1}^{*}}\right)\right] & = & \frac{1}{v_{1}^{*}}\left(1-\frac{v_{1}^{*}}{v(t)}\right)\left(N_{d}\delta\int_{0}^{h_{2}}f_{2}(\tau)y(t-\tau)d\tau-\mu v(t)\right).\end{eqnarray*}
Since, we have $\mu=\frac{N_{d}\delta y_{1}^{*}}{v_{1}^{*}}$, it
follows \begin{align}
\frac{d}{dt}\left[g\left(\frac{v(t)}{v_{1}^{*}}\right)\right] & =\frac{1}{v_{1}^{*}}\left(1-\frac{v_{1}^{*}}{v(t)}\right)\int_{0}^{h_{2}}f_{2}(\tau)\left(N_{d}\delta y(t-\tau)-N_{d}\delta y_{1}^{*}\frac{v(t)}{v_{1}^{*}}\right)d\tau\nonumber \\
 & =\frac{N_{d}\delta y_{1}^{*}}{v_{1}^{*}}\left(1-\frac{v_{1}^{*}}{v(t)}\right)\int_{0}^{h_{2}}f_{2}(\tau)\left(\frac{y(t-\tau)}{y_{1}^{*}}-\frac{v(t)}{v_{1}^{*}}\right)d\tau\nonumber \\
 & =\frac{N_{d}\delta y_{1}^{*}}{v_{1}^{*}}\left(1-\frac{1}{\tilde{v}_{t}}\right)\int_{0}^{h_{2}}f_{2}(\tau)\left(\tilde{y}_{t,\tau}-\tilde{v}_{t}\right)d\tau\nonumber \\
 & =\frac{N_{d}\delta y_{1}^{*}}{v_{1}^{*}}\int_{0}^{h_{2}}f_{2}(\tau)\left(\tilde{y}_{t,\tau}-\tilde{v}_{t}-\frac{\tilde{y}_{t,\tau}}{\tilde{v}_{t}}+1\right)d\tau.\label{eq:v2}\end{align}
Now, we see\begin{align}
 & \frac{d\overline{U}_{1}(t)}{dt}\nonumber \\
 & ={\displaystyle \int_{0}^{h_{1}}f_{1}(\tau)\left[g\left({\displaystyle \frac{x(t)v(t)}{x_{1}^{*}v_{1}^{*}}}\right)-g\left({\displaystyle \frac{x(t-\tau)v(t-\tau)}{x_{1}^{*}v_{1}^{*}}}\right)\right]d\tau}+{\displaystyle \int_{0}^{h_{2}}f_{2}(\tau)\left[g\left(\frac{y(t)}{y_{1}^{*}}\right)-g\left(\frac{y(t-\tau)}{y_{1}^{*}}\right)\right]d\tau}\nonumber \\
 & =\int_{0}^{h_{1}}f_{1}(\tau)\left({\displaystyle \tilde{x}_{t}}\tilde{v}_{t}-\ln\left({\displaystyle \tilde{x}_{t}}\tilde{v}_{t}\right)-{\displaystyle \tilde{x}_{t,\tau}}\tilde{v}_{t,\tau}+\ln\left({\displaystyle \tilde{x}_{t,\tau}}\tilde{v}_{t,\tau}\right)\right)d\tau+{\displaystyle \int_{0}^{h_{2}}f_{2}(\tau)\left(\tilde{y}_{t}-\ln\tilde{y}_{t}-\tilde{y}_{t,\tau}+\ln\tilde{y}_{t,\tau}\right)d\tau}.\label{eq:ca2}\end{align}
Consequently, by adding (\ref{eq:x2-1})-(\ref{eq:ca2}), we obtain
\begin{align}
\frac{d}{dt}U_{1}(t) & =-\frac{dx(t)}{k_{d}v_{1}^{*}x_{1}^{*}}\left(1-\frac{x_{1}^{*}}{x(t)}\right)^{2}+C_{1}(t,\tau)+C_{2}(t),\label{eq:last21}\end{align}
where\begin{align}
C_{1}(t,\tau) & =\left(1-\tilde{x}_{t}\tilde{v}_{t}-\frac{1}{\tilde{x}_{t}}+\tilde{v}_{t}\right)\nonumber \\
 & +\int_{0}^{h_{1}}f_{1}(\tau)\left(\tilde{x}_{t,\tau}\tilde{v}_{t,\tau}-\frac{\tilde{x}_{t,\tau}\tilde{v}_{t,\tau}}{\tilde{y}_{t}}-\tilde{y}_{t}+1\right)d\tau+\int_{0}^{h_{2}}f_{2}(\tau)\left(\tilde{y}_{t,\tau}-\tilde{v}_{t}-\frac{\tilde{y}_{t,\tau}}{\tilde{v}_{t}}+1\right)d\tau\nonumber \\
 & +{\displaystyle \int_{0}^{h_{1}}f_{1}(\tau)\left({\displaystyle \tilde{x}_{t}}\tilde{v}_{t}-\ln\left({\displaystyle \tilde{x}_{t}}\tilde{v}_{t}\right)-{\displaystyle \tilde{x}_{t,\tau}}\tilde{v}_{t,\tau}+\ln\left({\displaystyle \tilde{x}_{t,\tau}}\tilde{v}_{t,\tau}\right)\right)d\tau}\nonumber \\
 & +{\displaystyle \int_{0}^{h_{2}}f_{2}(\tau)\left(\tilde{y}_{t}-\ln\tilde{y}_{t}-\tilde{y}_{t,\tau}+\ln\tilde{y}_{t,\tau}\right)d\tau},\label{eq:C_11}\end{align}
and \begin{align}
C_{2}(t) & =-\frac{1}{k_{d}x_{1}^{*}v_{1}^{*}}\left(py(t)z(t)-py_{1}^{*}z(t)\right)+\frac{p}{k_{d}x_{1}^{*}v_{1}^{*}q}\left(\frac{d}{dt}z(t)\right)\nonumber \\
 & =-\frac{1}{k_{d}x_{1}^{*}v_{1}^{*}}\left(py(t)z(t)-py_{1}^{*}z(t)\right)+\frac{1}{k_{d}x_{1}^{*}v_{1}^{*}}\left(py(t)z(t)-\frac{p}{q}bz(t)\right)\nonumber \\
 & =\frac{1}{k_{d}x_{1}^{*}v_{1}^{*}}pz(t)\left(y_{1}^{*}-\frac{b}{q}\right).\label{eq:C_12}\end{align}
Now we claim $C_{2}(t)\leq0$ for all $t>0$. Since we have $R_{1}\leq1$,
$s\leq d\frac{\mu}{k_{d}N_{d}}+\frac{k}{k_{d}}\delta\frac{b}{q}$
holds from (\ref{eq:brn}). Then\begin{equation}
\frac{k_{d}}{k\delta}\left(s-d\frac{\mu}{k_{d}N_{d}}\right)=y_{1}^{*}\leq\frac{b}{q},\label{eq:R_1leq1}\end{equation}
from (\ref{eq:EE1}) and hence, \begin{equation}
C_{2}(t)\leq0.\label{eq:last22}\end{equation}
For $C_{1}(t,\tau)$, it holds that\begin{align}
C_{1}(t,\tau) & =\left(1-\frac{1}{\tilde{x}_{t}}\right)+\int_{0}^{h_{1}}f_{1}(\tau)\left(-\frac{\tilde{x}_{t,\tau}\tilde{v}_{t,\tau}}{\tilde{y}_{t}}+1\right)d\tau+\int_{0}^{h_{2}}f_{2}(\tau)\left(-\frac{\tilde{y}_{t,\tau}}{\tilde{v}_{t}}+1\right)d\tau\nonumber \\
 & \qquad\qquad+\int_{0}^{h_{1}}f_{1}(\tau)\left(-\ln\left({\displaystyle \tilde{x}_{t}}\tilde{v}_{t}\right)+\ln\left({\displaystyle \tilde{x}_{t,\tau}}\tilde{v}_{t,\tau}\right)\right)d\tau+{\displaystyle \int_{0}^{h_{2}}f_{2}(\tau)\left(-\ln\tilde{y}_{t}+\ln\tilde{y}_{t,\tau}\right)d\tau}\nonumber \\
 & =\left(1-\frac{1}{\tilde{x}_{t}}\right)+\int_{0}^{h_{1}}f_{1}(\tau)\left(-\frac{\tilde{x}_{t,\tau}\tilde{v}_{t,\tau}}{\tilde{y}_{t}}+1\right)d\tau+\int_{0}^{h_{2}}f_{2}(\tau)\left(-\frac{\tilde{y}_{t,\tau}}{\tilde{v}_{t}}+1\right)d\tau\nonumber \\
 & \qquad\qquad+\int_{0}^{h_{1}}f_{1}(\tau)\left(\ln\frac{1}{{\displaystyle \tilde{x}_{t}}}+\ln\frac{{\displaystyle \tilde{x}_{t,\tau}}\tilde{v}_{t,\tau}}{\tilde{y}_{t}}\right)d\tau+{\displaystyle \int_{0}^{h_{2}}f_{2}(\tau)\left(\ln\frac{\tilde{y}_{t,\tau}}{\tilde{v}_{t}}\right)d\tau}\nonumber \\
 & =\left(1-\frac{1}{\tilde{x}_{t}}+\ln\frac{1}{{\displaystyle \tilde{x}_{t}}}\right)\nonumber \\
 & +\int_{0}^{h_{1}}f_{1}(\tau)\left(-\frac{\tilde{x}_{t,\tau}\tilde{v}_{t,\tau}}{\tilde{y}_{t}}+1+\ln\frac{{\displaystyle \tilde{x}_{t,\tau}}\tilde{v}_{t,\tau}}{\tilde{y}_{t}}\right)d\tau+\int_{0}^{h_{2}}f_{2}(\tau)\left(-\frac{\tilde{y}_{t,\tau}}{\tilde{v}_{t}}+1+\ln\frac{\tilde{y}_{t,\tau}}{\tilde{v}_{t}}\right)d\tau\nonumber \\
 & =-g\left(\frac{1}{\tilde{x}_{t}}\right)-\int_{0}^{h_{1}}f_{1}(\tau)g\left(\frac{\tilde{x}_{t,\tau}\tilde{v}_{t,\tau}}{\tilde{y}_{t}}\right)d\tau-\int_{0}^{h_{2}}f_{2}(\tau)g\left(\frac{\tilde{y}_{t,\tau}}{\tilde{v}_{t}}\right)d\tau\leq0.\label{eq:last23}\end{align}

Consequently, $\frac{dU_{1}(t)}{dt}\leq0$ holds from (\ref{eq:last21}),
(\ref{eq:last22}) and (\ref{eq:last23}). Hence, every solution of
(\ref{eq:Mother-1}) tends to $M_{1}$, where $M_{1}$ is the largest
invariant subset in $\left\{ \frac{dU_{1}(t)}{dt}=0\right\} $ with
respect to (\ref{eq:Mother-1}). We show that $M_{1}$ consists of
only the equilibrium $E_{1}$. Let $\left(x(t),y(t),v(t),z(t)\right)$
be the solution with initial function in $M_{1}$, then, it holds
that \begin{equation}
x(t)=x_{1}^{*},\frac{x(t-\tau)v(t-\tau)}{x_{1}^{*}v_{1}^{*}}=\frac{y(t)}{y_{1}^{*}}\text{ for almost }\tau\in[0,h_{1}]\text{ and }\frac{y(t-\tau)}{y_{1}^{*}}=\frac{v(t)}{v_{1}^{*}}\text{ for almost }\tau\in[0,h_{2}].\label{eq:m1}\end{equation}
 From the invariance of $M_{1}$, we have $\frac{d}{dt}x(t)=0$ and
it then follows that $v(t)=v_{1}^{*}$ for any $t$ from the first
equation of (\ref{eq:Mother-1}). From (\ref{eq:m1}), we obtain $y(t)=y(t-\tau)=y_{1}^{*}$
for any $t$ and then, $z(t)=0$ follows from the second equation
of (\ref{eq:Mother-1}). Therefore, the infected equilibrium without
immune response $E_{1}$ is globally attractive. Since we have $\frac{dU_{1}(t)}{dt}\leq0$
and $U_{1}(t)\geq U_{1}(t)-\overline{U}_{1}(t)$, the infected equilibrium
without immune response $E_{1}$ is stable by Hale and Lunel \cite[Section 5, Corollary 3.1]{MR1243878}.
Hence, the infected equilibrium without immune response $E_{1}$ is
globally asymptotically stable for $R_{1}\leq1<R_{0}$.

iii) We construct the following Lyapunov functional \begin{equation}
U_{2}(t)=\frac{1}{kv_{2}^{*}}g\left(\frac{x(t)}{x_{2}^{*}}\right)+\frac{y_{2}^{*}}{k_{d}x_{2}^{*}v_{2}^{*}}g\left(\frac{y(t)}{y_{2}^{*}}\right)+\frac{v_{2}^{*}}{N_{d}\delta y_{2}^{*}}g\left(\frac{v(t)}{v_{2}^{*}}\right)+\frac{pz_{2}^{*}}{k_{d}x_{2}^{*}v_{2}^{*}q}g\left(\frac{z(t)}{z_{2}^{*}}\right)+\overline{U}_{2}(t),\label{eq: Lyapunov2}\end{equation}
 where\[
\overline{U}_{2}(t)=\int_{0}^{h_{1}}f_{1}(\tau)\int_{t-\tau}^{t}g\left(\frac{x(s)v(s)}{x_{2}^{*}v_{2}^{*}}\right)dsd\tau+\int_{0}^{h_{2}}f_{2}(\tau)\int_{t-\tau}^{t}g\left(\frac{y(s)}{y_{2}^{*}}\right)dsd\tau.\]
Similar to (\ref{eq:x2-1}), we obtain \begin{equation}
\frac{d}{dt}\left[g\left(\frac{x(t)}{x_{2}^{*}}\right)\right]=-\frac{dx(t)}{x_{2}^{*}}\left(1-\frac{x_{2}^{*}}{x(t)}\right)^{2}+kv_{2}^{*}\left(1-\overline{x}_{t}\overline{v}_{t}-\frac{1}{\overline{x}_{t}}+\overline{v}_{t}\right).\label{eq:x3}\end{equation}
We also obtain \begin{align*}
\frac{d}{dt}\left[g\left(\frac{y(t)}{y_{2}^{*}}\right)\right] & =\frac{1}{y_{2}^{*}}\left(1-\frac{y_{2}^{*}}{y(t)}\right)\left(k_{d}\int_{0}^{h_{1}}f_{1}(\tau)x(t-\tau)v(t-\tau)d\tau-\delta y(t)-py(t)z(t)\right).\end{align*}
Since we have $\delta y_{2}^{*}=k_{d}x_{2}^{*}v_{2}^{*}-py_{2}^{*}z_{2}^{*}$,
it holds \[
\delta=\frac{1}{y_{2}^{*}}\left(k_{d}x_{2}^{*}v_{2}^{*}-py_{2}^{*}z_{2}^{*}\right).\]
Then\begin{align}
 & \frac{d}{dt}\left[g\left(\frac{y(t)}{y_{2}^{*}}\right)\right]\nonumber \\
 & =\frac{1}{y_{2}^{*}}\left(1-\frac{y_{2}^{*}}{y(t)}\right)\left(k_{d}\int_{0}^{h_{1}}f_{1}(\tau)x(t-\tau)v(t-\tau)d\tau-\frac{1}{y_{2}^{*}}\left(k_{d}x_{2}^{*}v_{2}^{*}-py_{2}^{*}z_{2}^{*}\right)y(t)-py(t)z(t)\right)\nonumber \\
 & =\frac{1}{y_{2}^{*}}\left(1-\frac{y_{2}^{*}}{y(t)}\right)\left[\int_{0}^{h_{1}}f_{1}(\tau)\left(k_{d}x(t-\tau)v(t-\tau)-k_{d}x_{2}^{*}v_{2}^{*}\frac{y(t)}{y_{2}^{*}}\right)d\tau+\left(pz_{2}^{*}y(t)-py(t)z(t)\right)\right]\nonumber \\
 & =\frac{1}{y_{2}^{*}}\left(1-\frac{y_{2}^{*}}{y(t)}\right)\int_{0}^{h_{1}}f_{1}(\tau)\left(k_{d}x(t-\tau)v(t-\tau)-k_{d}x_{2}^{*}v_{2}^{*}\frac{y(t)}{y_{2}^{*}}\right)d\tau\nonumber \\
 & +\frac{1}{y_{2}^{*}}\left(1-\frac{y_{2}^{*}}{y(t)}\right)\left(pz_{2}^{*}y(t)-py(t)z(t)\right)\nonumber \\
 & =\frac{k_{d}x_{2}^{*}v_{2}^{*}}{y_{2}^{*}}\left(1-\frac{y_{2}^{*}}{y(t)}\right)\int_{0}^{h_{1}}f_{1}(\tau)\left(\frac{x(t-\tau)v(t-\tau)}{x_{2}^{*}v_{2}^{*}}-\frac{y(t)}{y_{2}^{*}}\right)d\tau\nonumber \\
 & +pz_{2}^{*}\left(1-\frac{y_{2}^{*}}{y(t)}\right)\left(\frac{y(t)}{y_{2}^{*}}-\frac{y(t)z(t)}{y_{2}^{*}z_{2}^{*}}\right)\nonumber \\
 & =\frac{k_{d}x_{2}^{*}v_{2}^{*}}{y_{2}^{*}}\left(1-\frac{1}{\overline{y}_{t}}\right)\int_{0}^{h_{1}}f_{1}(\tau)\left(\overline{x}_{t,\tau}\overline{v}_{t,\tau}-\overline{y}_{t}\right)d\tau+pz_{2}^{*}\left(1-\frac{1}{\overline{y}_{t}}\right)\left(\overline{y}_{t}-\overline{y}_{t}\overline{z}_{t}\right)\nonumber \\
 & =\frac{k_{d}x_{2}^{*}v_{2}^{*}}{y_{2}^{*}}\int_{0}^{h_{1}}f_{1}(\tau)\left(\overline{x}_{t,\tau}\overline{v}_{t,\tau}-\overline{y}_{t}-\frac{\overline{x}_{t,\tau}\overline{v}_{t,\tau}}{\overline{y}_{t}}+1\right)d\tau+pz_{2}^{*}\left(\overline{y}_{t}-1\right)\left(1-\overline{z}_{t}\right).\label{eq:y3}\end{align}
Similar to (\ref{eq:v2}), we also obtain\begin{equation}
\frac{d}{dt}\left[g\left(\frac{v(t)}{v_{2}^{*}}\right)\right]=\frac{N_{d}\delta y_{2}^{*}}{v_{2}^{*}}\int_{0}^{h_{2}}f_{2}(\tau)\left(\overline{y}_{t,\tau}-\overline{v}_{t}-\frac{\overline{y}_{t,\tau}}{\overline{v}_{t}}+1\right)d\tau.\label{eq:v3}\end{equation}
Let us calculate 

\begin{align}
\frac{d}{dt}\left[g\left(\frac{z(t)}{z_{2}^{*}}\right)\right] & =\frac{1}{z_{2}^{*}}\left(1-\frac{z_{2}^{*}}{z(t)}\right)\left(qy(t)z(t)-bz(t)\right)=\frac{1}{z_{2}^{*}}\left(1-\frac{z_{2}^{*}}{z(t)}\right)\left(qy(t)z(t)-qy_{2}^{*}z(t)\right).\nonumber \\
 & =\left(\frac{z(t)}{z_{2}^{*}}-1\right)\left(qy(t)-qy_{2}^{*}\right)\nonumber \\
 & =qy_{2}^{*}\left(\frac{z(t)}{z_{2}^{*}}-1\right)\left(\frac{y(t)}{y_{2}^{*}}-1\right)\nonumber \\
 & =qy_{2}^{*}\left(\overline{z}_{t}-1\right)\left(\overline{y}_{t}-1\right).\label{eq:z3}\end{align}
Similar to (\ref{eq:ca2}), we obtain\begin{align}
 & \frac{d\overline{U}_{2}(t)}{dt}\nonumber \\
 & ={\displaystyle \int_{0}^{h_{1}}f_{1}(\tau)\left(\overline{x}_{t}\overline{v}_{t}-\ln\left(\overline{x}_{t}\overline{v}_{t}\right)-\overline{x}_{t,\tau}\overline{v}_{t,\tau}+\ln\left(\overline{x}_{t,\tau}\overline{v}_{t,\tau}\right)\right)d\tau}+{\displaystyle \int_{0}^{h_{2}}f_{2}(\tau)\left(\overline{y}_{t}-\ln\overline{y}_{t}-\overline{y}_{t,\tau}+\ln\overline{y}_{t,\tau}\right)d\tau}.\label{eq:ca3}\end{align}

Consequently, by adding (\ref{eq:x3})-(\ref{eq:ca3}), we obtain
\begin{align}
\frac{d}{dt}U_{2}(t) & =-\frac{dx(t)}{x_{2}^{*}}\left(1-\frac{x_{2}^{*}}{x(t)}\right)^{2}+C_{3}(t,\tau),\label{eq:last31}\end{align}
where \begin{align*}
 & C_{3}(t,\tau)\\
 & =\left(1-\overline{x}_{t}\overline{v}_{t}-\frac{1}{\overline{x}_{t}}+\overline{v}_{t}\right)\\
 & +\int_{0}^{h_{1}}f_{1}(\tau)\left(\overline{x}_{t,\tau}\overline{v}_{t,\tau}-\overline{y}_{t}-\frac{\overline{x}_{t,\tau}\overline{v}_{t,\tau}}{\overline{y}_{t}}+1\right)d\tau+\int_{0}^{h_{2}}f_{2}(\tau)\left(\overline{y}_{t,\tau}-\overline{v}_{t}-\frac{\overline{y}_{t,\tau}}{\overline{v}_{t}}+1\right)d\tau\\
 & +\left[\frac{y_{2}^{*}}{k_{d}x_{2}^{*}v_{2}^{*}}pz_{2}^{*}\left(\overline{y}_{t}-1\right)\left(1-\overline{z}_{t}\right)+\frac{pz_{2}^{*}}{k_{d}x_{2}^{*}v_{2}^{*}q}qy_{2}^{*}\left(\overline{z}_{t}-1\right)\left(\overline{y}_{t}-1\right)\right]\\
 & +{\displaystyle \int_{0}^{h_{1}}f_{1}(\tau)\left(\overline{x}_{t}\overline{v}_{t}-\ln\left(\overline{x}_{t}\overline{v}_{t}\right)-\overline{x}_{t,\tau}\overline{v}_{t,\tau}+\ln\left(\overline{x}_{t,\tau}\overline{v}_{t,\tau}\right)\right)d\tau}+{\displaystyle \int_{0}^{h_{2}}f_{2}(\tau)\left(\overline{y}_{t}-\ln\overline{y}_{t}-\overline{y}_{t,\tau}+\ln\overline{y}_{t,\tau}\right)d\tau}\\
 & =\left(1-\overline{x}_{t}\overline{v}_{t}-\frac{1}{\overline{x}_{t}}+\overline{v}_{t}\right)\\
 & +\int_{0}^{h_{1}}f_{1}(\tau)\left(\overline{x}_{t,\tau}\overline{v}_{t,\tau}-\overline{y}_{t}-\frac{\overline{x}_{t,\tau}\overline{v}_{t,\tau}}{\overline{y}_{t}}+1\right)d\tau+\int_{0}^{h_{2}}f_{2}(\tau)\left(\overline{y}_{t,\tau}-\overline{v}_{t}-\frac{\overline{y}_{t,\tau}}{\overline{v}_{t}}+1\right)d\tau\\
 & +{\displaystyle \int_{0}^{h_{1}}f_{1}(\tau)\left(\overline{x}_{t}\overline{v}_{t}-\ln\left(\overline{x}_{t}\overline{v}_{t}\right)-\overline{x}_{t,\tau}\overline{v}_{t,\tau}+\ln\left(\overline{x}_{t,\tau}\overline{v}_{t,\tau}\right)\right)d\tau}+{\displaystyle \int_{0}^{h_{2}}f_{2}(\tau)\left(\overline{y}_{t}-\ln\overline{y}_{t}-\overline{y}_{t,\tau}+\ln\overline{y}_{t,\tau}\right)d\tau}.\end{align*}
Similar to (\ref{eq:C_11}), we see \begin{equation}
C_{3}(t,\tau)=-g\left(\frac{1}{\overline{x}_{t}}\right)-\int_{0}^{h_{1}}f_{1}(\tau)g\left(\frac{\overline{x}_{t,\tau}\overline{v}_{t,\tau}}{\overline{y}_{t}}\right)d\tau-\int_{0}^{h_{2}}f_{2}(\tau)g\left(\frac{\overline{y}_{t,\tau}}{\overline{v}_{t}}\right)d\tau\leq0.\label{eq:last32}\end{equation}

Thus, $\frac{dU_{2}(t)}{dt}\leq0$ holds from (\ref{eq:last31}) and
(\ref{eq:last32}). Hence, the solution of system (\ref{eq:Mother-1})
limit to $M_{2}$, where $M_{2}$ is the largest invariant subset
in $\left\{ \frac{dU_{2}(t)}{dt}=0\right\} $ with respect to (\ref{eq:Mother-1}).
We show that $M_{2}$ consists of only the equilibrium $E_{2}$. Let
$\left(x(t),y(t),v(t),z(t)\right)$ be the solution with initial function
in $M_{2}$, then it holds that \begin{equation}
x(t)=x_{2}^{*},\frac{x(t-\tau)v(t-\tau)}{x_{2}^{*}v_{2}^{*}}=\frac{y(t)}{y_{2}^{*}}\text{ for almost }\tau\in[0,h_{1}]\text{ and }\frac{y(t-\tau)}{y_{2}^{*}}=\frac{v(t)}{v_{2}^{*}}\text{ for almost }\tau\in[0,h_{2}].\label{eq:m2}\end{equation}
From the invariance of $M_{2}$, we have $\frac{d}{dt}x(t)=0$ and
it then follows that $v(t)=v_{2}^{*}$ for any $t$ from the first
equation of (\ref{eq:Mother-1}). From (\ref{eq:m2}), we obtain $y(t)=y(t-\tau)=y_{2}^{*}$
for any $t$ and then, $z(t)=z_{2}^{*}$ follows from the second equation
of (\ref{eq:Mother-1}). Therefore, the infected equilibrium with
immune response $E_{2}$ is globally attractive. Since we have $\frac{dU_{2}(t)}{dt}\leq0$
and $U_{2}(t)\geq U_{2}(t)-\overline{U}_{2}(t)$, the infected equilibrium
with immune response $E_{2}$ is stable by Hale and Lunel \cite[Section 5, Corollary 3.1]{MR1243878}.
Hence, the infected equilibrium with immune response $E_{2}$ is globally
asymptotically stable for $R_{1}>1$.

Finally, the proof of this theorem is complete. \qed 
\end{pf}

\section{Applications}

Our approach is applicable for discrete delay models. Zhu and Zou
\cite{HZhu2009} studied the following viral infection model with
cell mediated immunity.\begin{equation}
\begin{cases}
\frac{d}{dt}x(t)=s-dx(t)-kx(t)v(t),\\
\frac{d}{dt}y(t)=k\textrm{e}^{-\delta\tau}x(t-\tau)v(t-\tau)-\delta y(t)-py(t)z(t),\\
\frac{d}{dt}v(t)=N\delta y(t)-\mu v(t),\\
\frac{d}{dt}z(t)=qy(t)z(t)-bz(t),\end{cases}\label{eq:zhuzou}\end{equation}
with the initial conditions $x(\theta)=\varphi_{1}(\theta),y(0)=y_{0},v(\theta)=\varphi_{3}(\theta),z(0)=z_{0}\text{ for }\theta\in[-\tau,0],$
where $\varphi_{i}(\theta)\in C([-\tau,0],\mathbb{R}_{+}),i=1,3,y_{0}\geq0$
and $z_{0}\geq0$. All parameters are positive constant.

For (\ref{eq:zhuzou}), similar to (\ref{eq:Mother-1}), there exist
three possible equilibria. From (\ref{eq:brn}) and (\ref{eq:brn2}),
the basic reproduction number for viral infection and for CTL response
are given by \[
\overline{R}_{0}=\frac{s}{d\frac{\mu}{k\textrm{e}^{-\delta\tau}N}}\text{ and }\overline{R}_{1}=\frac{s}{d\frac{\mu}{k\textrm{e}^{-\delta\tau}N_{d}}+\textrm{e}^{\delta\tau}\delta\frac{b}{q}},\]
respectively. There exist the uninfected equilibrium $E_{0}(x_{0},0,0,0),x_{0}=\frac{s}{d}$,
the infected equilibrium without immune response $\overline{E}_{1}(\overline{x}_{1}\overline{y}_{1},\overline{v}_{1},0)$
($\overline{x}_{1},\overline{y}_{1},\overline{v}_{1}>0$) if $\overline{R}_{0}>1$
and the infected equilibrium with immune response $\overline{E}_{2}(\overline{x}_{2}\overline{y}_{2},\overline{v}_{2},\overline{z}_{2})$
($\overline{x}_{2},\overline{y}_{2},\overline{v}_{2},\overline{z}_{2}>0$)
if $\overline{R}_{1}>1$ (see also \cite[Section 3]{HZhu2009}). 

Zhu and Zou \cite{HZhu2009} established the global asymptotic stability
of the uninfected equilibrium $E_{0}$ for $\overline{R}_{0}<1$.
Moreover, they obtained sufficient conditions for the local asymptotic
stability of infected equilibria $\overline{E}_{1}$ and $\overline{E}_{2}$
by analysis of associated characteristic equations. Complete global
dynamics for (\ref{eq:zhuzou}) is not clear and an open problem.
However, similar to Theorem \ref{thm:Main} in Section 3, we establish
the following result.
\begin{thm}
\label{thm:app}i) If $\overline{R}_{0}\leq1$, then the uninfected
equilibrium $E_{0}$ for (\ref{eq:zhuzou}) is globally asymptotically
stable.

ii) Assume $y_{0}+\int_{0}^{h_{1}}f_{1}(\tau)\varphi_{1}(-\tau)\varphi_{3}(-\tau)d\tau>0$,
or $\varphi_{3}(0)>0$. If $\overline{R}_{1}\leq1<\overline{R}_{0}$,
then the infected equilibrium without immune response $\overline{E}_{1}$
for (\ref{eq:zhuzou}) is globally asymptotically stable.

iii) Assume $z_{0}>0$ and either $y_{0}+\int_{0}^{h_{1}}f_{1}(\tau)\varphi_{1}(-\tau)\varphi_{3}(-\tau)d\tau>0$,
or $\varphi_{3}(0)>0$. If $\overline{R}_{1}>1$, then the infected
equilibrium with immune response $\overline{E}_{2}$ for (\ref{eq:zhuzou})
is globally asymptotically stable.
\end{thm}
Zhu and Zou \cite[Theorems 3.3, 3.4]{HZhu2009} showed that the infected
equilibrium without immune response $\overline{E}_{1}$ is locally
asymptotically stable for $\overline{R}_{1}<1<\overline{R}_{0}$ and
the infected equilibrium with immune response $\overline{E}_{2}$
is locally asymptotically stable for $\overline{R}_{1}>1$ if the
intracellular delay $\tau$ satisfies a condition (see \cite[Theorem 3.4]{HZhu2009}).
However, by Theorem \ref{thm:app}, we establish that $\overline{E}_{1}$
is not only locally asymptotically stable but also globally asymptotically
stable for $\overline{R}_{1}\leq1<\overline{R}_{0}$. Moreover, $\overline{E}_{2}$
is globally asymptotically stable, whenever it exists, that is, $\overline{R}_{1}>1$.

\section{Discussion}

In this paper, we study global dynamics of delay differential equations
for a virus-immune interaction in \textit{vivo}. Two distributed time
delays represent the time needed for infection of cell and virus replication.
Stability analysis for (\ref{eq:Mother-1}) with discrete intracellular
delay was carried out by Li and Shu \cite{Michael1} and Zhu and Zou
\cite{HZhu2009}. Li and Shu \cite{Michael1} studied a viral infection
model which ignores the immune response to the viral infection and
showed that their model always admits an equilibrium which is globally
asymptotically stable. Recently, Li and Shu \cite{li:2434} has investigated
a general viral infection model with distributed delay which also
does not incorporate the immune response. Zhu and Zou \cite{HZhu2009}
established global stability of an uninfected equilibrium and obtained
sufficient conditions for local asymptotic stability of two infected
equilibria when the distributed delay in (\ref{eq:Mother-1}) is given
by a discrete. Zhu and Zou \cite{HZhu2009} did not address the global
stability of two infected equilibria for their model.

To obtain an integrated view for the virus-immune interaction dynamics
in \textit{vivo}, we investigate the global stability of (\ref{eq:Mother-1})
by employing the method of Lyapunov functionals which are motivated
by McClusky \cite{McCluskey201055} for delayed epidemic models. (\ref{eq:Mother-1})
has three possible equilibria, an uninfected equilibrium and two infected
equilibria with or without immune response. A combination of the basic
reproduction number for viral infection $R_{0}$ and for CTL response
$R_{1}$, defined by (\ref{eq:brn}) and (\ref{eq:brn2}), respectively,
determine the existence of these equilibria. Moreover, they also fully
determine the global dynamics of the model. The uninfected equilibrium
$E_{0}$ is globally asymptotically stable if $R_{0}\leq1$ and the
viruses are cleared. The infected equilibrium without immune response
$E_{1}$ is globally asymptotically stable if $R_{1}\leq1<R_{0}$
and the infection becomes chronic but with no persistent immune response.
The infected equilibrium with immune response $E_{2}$ is globally
asymptotically stable if $R_{1}>1$ and the infection becomes chronic
with immune response. Theorem \ref{thm:Main} is an extension result
of the global stability results in Pr$\ddot{\textrm{u}}$ss et al.,
\cite{MR2460257} and Li and Shu \cite{Michael1}. Moreover, we improve
stability results in Zhu and Zou \cite{HZhu2009} (see Section 4).

We see that virus eventually persists if $R_{0}>1$, because the infected
equilibrium $E_{1}$ or $E_{2}$ is globally asymptotically stable
in this case. The infected equilibrium without immune response $E_{1}$
is globally asymptotically stable and the immune response does not
work for $R_{0}>1\geq R_{1}$. On the other hand, the immune response
is activated and there exist two infected equilibria $E_{1}=(x_{1}^{*},y_{1}^{*},v_{1}^{*},0)$
and $E_{2}=(x_{2}^{*},y_{2}^{*},v_{2}^{*},z_{2}^{*})$ for $R_{1}>1$.
Moreover, in these equilibria, one can see that the relations $x_{1}^{*}<x_{2}^{*}$
and $y_{2}^{*}<y_{1}^{*}$ hold due to the effect of immunity (see
also Remark \ref{rem:equ_cond}). Therefore, the global stability
of $E_{2}$ for $R_{1}>1$ indicates that the immune activation has
a positive role in the reduction of the infected cells and the increasing
of the uninfected cells for $R_{1}>1$.
\begin{ack}
The author thanks the referee for the valuable suggestions and numerous
comments which led to a significant improvement on the original manuscript.
The work of this paper was partially done during the author's visit
as members of International Research Training Group 1529 in TU Darmstadt
on February, 2010. Finally, the author is grateful to Professor Rico
Zacher for the constructive comments and for bringing the paper \cite{MR2460257}
to his attention at the international workshop \textquotedblleft{}Mathematical
Fluid Dynamics\textquotedblright{} held in Waseda University, Tokyo
on 8-16 March 2010.
\end{ack}

\end{document}